\documentclass[a4,12pt]{amsart}
\usepackage{amsmath, amsthm, amssymb, setspace, graphics}
\usepackage{mdframed}
\usepackage[all]{xy}
\usepackage{url}
 \doublespace

\newtheorem{thm}{Theorem}[subsection] 
\theoremstyle{definition} 
\theoremstyle{remark} 
\newtheorem{remark}[thm]{Remark}

\newcommand{\p}{\partial}

\begin{document}

\title{{Remark on symplectic relative Gromov-Witten invariants and degeneration formula }}
\author{An-Min Li}


\begin{abstract}
   In this note, we give  item-by-item responses to the criticisms raised in   \cite{TZ}  by Tehrani ad Zinger on our paper \cite{LR}.  We illuminate the main ideas and contributions in \cite{LR} in section 2,  itemize the responses to issues  raised in \cite{TZ} and conclude
  that  we have provided a complete proof of the degeneration formula in   our published  paper   \cite{LR} and its more detailed versions in arXiv.   In \cite{TZ},  the authors made an   effort in comparing the methods and ideas in \cite{LR} vs \cite{IP-1} \cite{IP-2}, but  their criticisms
 on   \cite{LR}  are based on their  own lack of sufficient understanding of \cite{LR}.
\end{abstract}

\maketitle
\tableofcontents

\section{Introduction}\label{sect_1}

More than fifteen years ago, Yongbin Ruan and I developed a theory of relative Gromov-Witten invariants and degeneration formula (see \cite{IP-1} and  \cite{IP-2} for a different approach and {readers may refer to Remark 2.3(\cite{TZ}) for a detailed documentation}).
Since then, our formula has been applied and tested many times, for example, an algebraic treatment of our theory was developed by Jun Li \cite{L}.
 Recently, there was an
article \cite{TZ} casting certain doubt of our results. As far as we understood, the author did not question the correctness of our result as well as the  effectiveness of our approach.
The dispute is if we provided enough detail which qualifies it as a complete proof. Incidently, we do have a much longer version of paper available online which is predated the published
version and contains much more details. For the published version, the referee forced us to cut off 40 pages of material which he considered to be standard. Namely, our
long version was considered to contain too much details. Ironically, fifteen years later the readers of different generation complains that our shorter published version has too few
details. When the article \cite{TZ} was first circulated in a large mailing list, we informed the author our long version. They refused to consider it! Since the issue of enough detail
is precisely the center of dispute, we question the author's fairness in treating our work.

Nevertheless, we feel that it is our responsibility to answer any questions anyone may have
for our work. This is the purpose of this article. Originally, we hope to result the difference through a private discussion. Unfortunately, we were not given such an opportunity.
We regret that we have to response to the article \cite{TZ}  in public. When the paper \cite{LR} was  in preparation, Yongbin Ruan was in the process of moving to other areas of mathematics. Since then, he has invented several important areas such as Chen-Ruan
cohomology and FJRW-theory. The theory of relative Gromov-Witten
theory and degeneration formula was written up by myself which I will take full responsibility for its correctness and completeness. Instead of dragging him back from his current important works, it is more appropriate for me to respond to all the criticisms in \cite{TZ}.

\subsection{Background}\label{sect_1.1}

  Our paper was written more than fifteen years ago for a different generation of mathematicians. Every paper assumes reader's  familiarity of certain basic or standard material. Our paper is
  no exception. To help current younger generation to understand our paper and its production, it is very important to discuss the background of our paper and what was considered to be
  standard material then. Our approach was an adaption of the so-called {\em neck stretching  technique}. This technique  has been  developed  in gauge theory in the late 80 and early 90's under the name
  of analysis on manifolds with cylindric end or { $L^2$-moduli space theory}. Floer homology is such an example. When our paper was prepared, this technique was already quite standard in the gauge theory community.  There are several  books on this technique,  for  example,   two books on $L^2$-moduli space theory by Mrowka-Morgan-Ruberman(\cite{MMR}) and by Taubes (\cite{Taubes}).
  Donaldson's beautiful  manuscript  on Floer homology was   available to the public (the actual book  \cite{D} was published in 2004) .
   Our paper was prepared between 1996 and 1998. We certainly assumed  readers'  familiarity in  some basic or standard knowledge of
  this neck stretching   technique. Nevertheless, {  in our first version of the paper, we did provide  a rather complete version on this
  technique in Gromov-Witten theory.(cf. \cite{LR} Version 1.)}

After our paper was submitted to Inventiones mathematicae, it went through  a long refereeing process  and  several
revisions were produced. One of disputes with referee is  exactly that we want to assume less  background on the   neck stretching  technique while he/she wants to assume more.
There is no way we can foresee that next generation of mathematicians does not know so much about this standard  technique.
We welcome any effort to rewrite
  our theory in the new language which young people are more familiar with it. However, the article \cite{TZ} is different. Namely, they challenged the completeness of our proof. However, we feel that it is the challengers'
  responsibility to actually  be familiar with  the  technique we applied  (a standard technique in  90's) and to understand  our proof before making the judgement.
{  In this response, we discuss this issue in details in \S\ref{sect_2.4} and \S\ref{sect_3.1}.}

{  Another crucial  technical issue in \cite{LR} involved is to define invariants using the virtual  neighbourhood technique. Unlike the $L^2$-moduli spaces in the Yang-Mills theory, which had been well estabilished by that  time, the developement of virtual techniques  for  Gromov-Witten theory was just at the early stage, and even after 15 years' today. That is why  SCGP (Stony Brook) has a half-year program on the foundation of Gromov-Witten theory. As we know, there had been several different approaches  at that  time, such as Fukaya-Ono(\cite{FO}), Li-Tian(\cite{LT}), Liu-Tian(\cite{LiuT}),
Ruan(\cite{R}), Siebert(\cite{S}), just to name a few.
In \cite{LR}, we provide a completely new approach to this issue: we show that the invariants can be defined via the integration on the {\em top} stratum  in the sense of virtual neighbourhood.  It turns out  that this new approach is  very effective   in many  later applications. Here again, we feel that the authors of \cite{TZ} are not  very familiar with this new viewpoint  and hence have no idea of the efficiency and hence the  correctness of our approach in \cite{LR}. We will  discuss this issue
in more details in \S\ref{sect_2.2}.}

The article \cite{TZ}  posed 16 specific questions. After studying their questions carefully,
  we concluded that our proof  is  COMPLETE. In fact,  in this note we shall demonstrate that most of
  their complaints and  criticisms   of  \cite{LR} are resulted precisely from the author's lack of basic
  understanding  of our approach.  It will certainly  take a while for the authors to really understand the paper \cite{LR}. We sincerely hope that when  they finally understand the main techniques in \cite{LR} and will be able to rewrite in their own language, they won't claim that they provide a complete proof of those theorems in \cite{LR}.

\section{Response}\label{sect_2}

     In this section, we  will answer explicitly the 16 questions  posed by the authors of  \cite{TZ}.
We remark that when we were preparing \cite{LR},  as a standard practice, we write  it as a research paper rather than a textbook. Hence, we illuminate    the proofs in \cite{LR}  that we believe  is enough for experts  to understand, not to provide a training wheel for all  those not  familiar with the knowledge accumulated over the last decade.   Any detail that can be
routinely filled would not be presented here.
Here we focus mainly on the new ideas in \cite{LR}.

We like to point out  that comments of Tehrani and Zinger (T \& Z for short), (LR1)-(LR16) (P15 in \cite{TZ}),
 may be classified as the following  three types.
\begin{enumerate}
\item   comments that  are related to standard materials to the subjects, e.g
(LR4), (LR6), (LR9), (LR11), (LR12);
\item comments  that are on some minor typo or overlooks that can be easily
fixed by diligent readers, e.g, (LR1), (LR2), (LR5),  (LR15);
\item  comments  that are on the mathematical techniques developed in \cite{LR}.  Reading \cite{TZ}, it is clear to us that T \& Z   either misunderstood or did not understand  at all of these mathematical techniques. They often  made   ridiculous comments on the mathematics
in \cite{LR},  even on  some of materials that  are already well known nowadays.  For example, it is clear that T \& Z are not familiar at all about the  Fredholm analysis  and the compact properties of  the $L^2$-moduli spaces when there are  certain Bott-Morse type equations involved (cf. (LR4) and  (LR5)).  Other  similar comments include (LR3),
LR(7), LR(8), (LR10), (LR13), (LR14) and (LR16).
\end{enumerate}
Since T \& Z  make  many incorrect comments  even on what we believe had already been quite standard materials by the time when the paper  \cite{LR}  was written, we can't help to question their  expertise on this  topic to judge the correctness of the paper \cite{LR}.
Moreover, on the main contributions of \cite{LR} to the development of the degenerate formulae and their application to symplectic topology and birational symplectic geometry,
we  remake  that

\begin{enumerate}
\item  T \& Z understand neither the approach in \cite{LR}  nor the essence of  proofs therein;
\item   they simple made their wishful and often ignorant   judgements based on their  self-claimed   righteous mathematical viewpoint.
 \end{enumerate}
Due to these, we will begin  with an outline of the approach  in \cite{LR} to the symplectic sum formula and highlight the new points of that paper in \S \ref{sect_2.2},\S\ref{sect_2.3} and \S\ref{sect_2.4}.  We let mathematics itself   in this note to speak for itself.

\subsection{Outline of approach to symplectic sum formula (\cite{LR})}\label{sect_2.1} We use the same notations as in \cite{LR}.
Let $(M, \omega)$ be a compact symplectic manifold of dimension $2n+2$, and $\widetilde M = H^{-1}(0)$ for a  local Hamiltonian function
as in the beginning of Section 3 in \cite{LR}.  Under the assumption that the  Hamiltonian vector field $X_H$
 generates a circle action on a neighborhood of $\widetilde M$, there
is a circle bundle  $\pi: \widetilde M \to Z =\widetilde M /S^1$ with a natural symplectic form $\tau_0$ on $Z$. We assume that $\widetilde M$ separates $M$ into two parts   to produce two cylindrical end  symplectic manifold  $M^+$ and $M^-$.  Collapsing the $S^1$-action
 at the infinity, we obtain the symplectic cuts $\overline{M}^+$ and  $\overline{M}^-$, both contain $Z$ as a  codimension two symplectic submanifold. We also consider the limiting manifold $M_\infty$ as we stretch the neck along $\widetilde M$.

To obtain and prove the symplectic sum formula, we began with the following strategies.
\begin{enumerate}
\item [(A)]  we relate the Gromov-Witten invariants  of $M$ with that of $M_\infty$(cf. Theorem 5.6 in \cite{LR}),
\item[(B)]  then we  relate the   Gromov-Witten invariants of $M_\infty$ with  the combination of relative invariants of $(\bar{M}^\pm,Z)$ (cf. Theorem 5.7 in \cite{LR}).
\end{enumerate}
Note that (A) and (B) yields the symplectic sum formula.
 For this purpose,
 \begin{enumerate}
 \item[(C)]  We introduce the relative moduli spaces for symplectic pairs $(\overline{M}^\pm,Z)$ (cf. Definition 3.14 in \cite{LR}) and  the moduli spaces on $M_\infty$ (cf. Definition 3.18  in \cite{LR});
 \item[(D)] Then we define the invariants for these moduli spaces, in particular,  including the relative GW invariants of  $(\overline{M}^\pm,Z)$.
\end{enumerate}
We will recall the main ideas to get (C) and (D)  in \S\ref{sect_2.3} and \S\ref{sect_2.2}.

\begin{remark}\label{remark_2.1.1}
{ We would like to mention our work on relative orbifold Gromov-Witten theory(\cite{CLSZ}). In \cite{CLSZ}, we  employ  a different approach to get  the symplectic sum formula for orbifold Gromov-Witten invariants. Instead, for a degeneration family of  symplectic orbifolds, we construct a degeneration family of moduli spaces of pseudo-holomorphic curves. In \cite{CLSZ} we then adapt  an integration argument to conclude the symplectic sum formula easily.
On  the moduli space level, T \& Z's approach in \cite{TZ} seems  very similar to the approach in \cite{CLSZ}.  }
\end{remark}

\subsection{The approach of defining invariants}\label{sect_2.2}
In \cite{LR}, we introduce a {\em new} approach for the  definition of Gromov-Witten invariants, which is not same as those virtual fundamental class approach in  the  existing literature.  By refinement of estimates of gluing maps,
we showed  that these  invariants can be defined via the {\em integration} on {\em top stratum} of the moduli space. As we know,  the  common way to define invariants is by using the  intersection theory. Of course,  one envisages  that one can define the invariants by considering  certain {\em virtual } intersection  theory for the  top stratum.  A dual approach is to define invariants via the integration, for example, see P. 227 in \cite{MS-2}.  However, in order to make sense of the integration theory, one needs to establish  the smoothness of compact moduli space.We know that the smoothness of the space  has only been achieved  very recently by several groups (including  my recent joint  work with Bohui Chen and Bai-Ling Wang).  However, in this paper, we shall  {\em avoid} the smoothness at lower strata, and prove that "the invariants defined by integration can be obtained by  integrating virtually on the  top stratum".
We believe that this  is a highly nontrivial and very useful statement (see the following Remark).

\begin{remark}\label{remark_2.2.1}
When   we  compare invariants mentioned in (A) and (B) ( in \S\ref{sect_2.1} ), we only need to compare the  top strata of moduli spaces in the virtual sense.
\end{remark}

\def \om{\overline{\mathcal M}}
\def \mc{\mathcal}
We now outline this key idea in \cite{LR}.
Let $\om$ be a compactified moduli space of the  top stratum $\mc M$. Set $\partial \mc M=\om\setminus \mc M$.
 For simplicity, we first assume the regularity holds for  $\om$
 and suppose that $\partial\mc M$  is of codimension $\geq 2$.  By gluing maps, we may give local coordinate charts for neighborhoods of any point $x\in \partial\mc M$.  At the bottom of  Page 204 in \cite{LR}, we explained  that  the integrand, which is {\em well defined} on
$\mc M$, behaviors well near $\partial M$ with respect to the local coordinate charts mentioned above. Hence, the invariants can be defined via the integrations on $\mc M$.   In fact, such a  strategy  is commonly used for singular spaces.  In order to
achieve this goal, we provided much more {\em refined} estimates for differential $\partial/\partial r$ of gluing maps,  for example, see  Lemma  4.7, 4.8, 4.9 and 4.13 in \cite{LR}. All these estimates were not appeared in any literature.   We would like to point out that in recent work of FOOO on the smoothness of  moduli space, they also consider the estimates of the similar type.

\begin{remark}\label{remark_2.2.2} We remark that if $\om$ is not regular, we should apply the virtual neighborhood technique, and then apply the above argument to the virtual neighborhood.
T \& Z commented ( LR16 in \cite{TZ})) that the above approach is not necessary since one can take the approach of intersections  etc.
We are shocked of  this kind of  naive viewpoint. This explained why  they either did not understand or  don't respect others' work though it is clear that T \& Z  has learnt a lot from \cite{LR}.
\end{remark}

\subsection{Relative moduli spaces and relative Gromov-Witten invariants}
\label{sect_2.3}
In \cite{LR},  we introduced the relative moduli space
of stable maps and define the relative GW invariants.

Let $(M,Z)$ be a symplectic pair as in \cite{LR} and let $\mc M$ be a relative moduli space.
By the time the paper was written, it is well known that  the compactification
argument from contact geometry would yield a codimension one boundary. This
issue is first resolved  {\em in \cite{LR}} by introducing the $\mathbb C^\ast$ action on the space of maps to rubber components.  {  For example, one might compare with \cite{IP-1} which is the first version of IP's series. In \cite{IP-1}, they did not acturally provide the precise construction of relative moduli spaces. However, from their statements, their invariants are only well defined up to chambers. This is
exactly  due to the boundary of  codimension one issue.}

Once we have this key observation, it remains to build up the moduli space $\om$
following the  following standard steps:
\begin{enumerate}
\item[(E)] compactification (\S3(\cite{LR})): we adapt the standard $L^2$-moduli space theory that
has been intensively developed for Chern-Simons theory, for example, we follow closely
with Donaldson' book (\cite{D}).  More details will be discussed in \S\ref{sect_2.4} and \S\ref{sect_3.1};
\item[(F)] regularization(\S 4.1(\cite{LR})): we use Ruan's argument to build up a global regularization;
\item[(G)] gluing theory(\S4.1(\cite{LR})): the standard package of gluing theory consists of the following steps
\begin{enumerate}
\item[(G1)]  the construction of gluing maps,
\item[(G2)] the injectivity and surjectivity of the gluing maps.
\end{enumerate}

In \cite{LR}, we explained  (G1) mainly for one {\em relative} nodal case, and we think  that  the generalisation to  lower strata are already standard. (cf. \S\ref{sect_3.2.3}). With understood,  we focused instead  on
\begin{enumerate}
\item[(G3)]  the  refined  and new estimates for  gluing maps, which was of  importance to our approach of defining the relative  invariants (cf. (H) below and \S\ref{sect_2.2});
\end{enumerate}

\item[(H)] to define  the invariants(\S4.2(\cite{LR})): this already explained in \S\ref{sect_2.2}; and again we emphasize that  this depends heavily on our refined and new  estimates for gluing maps (G3).
\end{enumerate}

\subsection{$L^2$ moduli space theory and Bott-Morse type Morse theory}
\label{sect_2.4}
The $L^2$ moduli space theory is perfectly suitable for  moduli spaces modeled on cylinder end domain.  It relates to study gradient flows of certain infinite dimensional Bott-Morse   function.  In
our case,  this  function is the smooth functional  $\mc A$ on  the smooth Banach manifold $W_r^2(S^1,\tilde M)$ (cf. \S3.1in \cite{LR}  for the definition of $\mc A$).
\begin{remark}\label{remark_2.4.1}
In the first version of the paper (see [arXiv 9803036v1],P52), $\mc A$ is defined up to
constants $\tau_0(T^2)$. Due to the opinion of referees, we cut the paper  into this new shorter version, in this version we only use a local version of
$\mc A$. In fact, mathematically, the global one is more convenient for applications.

\end{remark}
Using  the standard arguments, the crucial  part of the $L^2$ moduli space theory is to develop the Morse theory of $\mc A$ in the sense of Floer.
In \cite{HWZ}, authors studied the case that the critical points are isolated. While in our case, we show that $\mc A$ is of Bott-Morse type.  Hence the main goal of \S3.1 in \cite{LR}  is to generalize the theory of HWZ to the case of Bott-Morse type. The key ingredient  is Proposition 3.4 in \cite{LR}  which enabled  us  to prove Theorem 3.7 in \cite{LR}.

\begin{remark}\label{remark_2.4.2}
T \& Z questioned  about  Proposition 3.4 and Theorem 3.7 in \cite{LR}.  All of their questions   can be solved by   standard  infinite dimensional analysis for the  Bott-Morse type functional.
  We will explain in more details in \S\ref{sect_3.1}.
\end{remark}

\def \v{\vskip 0.1in}

\section{Response to TZ's comments}\label{sect_3}

In this section, we respond to  the comments (LR1-16) of T \& Z on \cite{LR}  according the topics discussed  in the the previous section.

\subsection{Compactification}\label{sect_3.1} We comment on (LR4-8) in this subsection.

As explained in \S\ref{sect_2.4}, the compactification of the relative moduli space  depends heavily on the
study of gradient flow equation  of  the functional $\mc A$, which is of Bott-Morse type.  T \& Z questioned  this soundness of this  technical part, in particular,
their comments (LR4) and (LR5) directed at  Proposition 3.4 and Theorem 3.7 respectively.  We suspect  that T\& Z are not familiar with  those  standard analysis involved with the study of  $L^2$ moduli spaces.

\subsubsection{}\label{sect_3.1.1}
On (LR4):  T\& Z complained that (1) we use the $L^2$ inner product to define a Riemann metric on the Banach manifold $W^2_r(S^1,SV)$; and  (2) the  infinite dimensional version of Morse lemma is used.
\v
{\bf Response:}
\v
(1) we have no idea what is wrong with $L^2$ inner product. In fact, this is {\em widely} used, and, for example, cf. the equation in   Page 32 of  \cite{D}, or Lemma 2.1.1 in \cite{MMR}.
\v
(2) the argument using Morse Lemma is also standard.  For example, one may be referred to Page 29 in \cite{D} to how this is applied for Yang-Mills Floer homology.   We thought and still think that  this part of the infinite dimensional  Morse theory is well known. There is no
need to explain further here about T\& Z's complaints.
\subsubsection{(LR5)}\label{sect_3.1.2}
On (LR5): T \& Z commented  that Theorem 3.7 in \cite{LR} (1) incorrectly stated;
(2) the proof including circular reasoning;  and  (3) the theorem can be justified in a few lines. They further claimed that
the statement after the proof of [LR, Theorem 3.7] does not make sense, because the constants there depend on the map $C\rightarrow \mathbb  R\times \tilde{M}$. In Proposition 5.8 of \cite{TZ} the constant is  $C_u$  which is depending  on $u$.
\v

{\bf Response:}
\v
(1)   Lemma 3.7 is a local theorem concerning with the behaviours of $u$ near $\infty$.  We note that $\mathbb C$ in the statement of  this Lemma  should be
$\mathbb C\setminus \mathbb D$. From the context of our paper it is easy to see that this is clearly a  editorial typo (in fact, we had
$\mathbb C\setminus\mathbb D$  in the first version of our paper [arXiv 9803036v1], Page 29  Theorem 3.7 in \cite{LR} ).
\v
 (2) We think that T\& Z totally misunderstood or could  not follow our proof: (cf. Remark 5.7 in \cite{TZ} ) T\& Z said  that we presupposed  "the flow stays within a small neighborhood $O_{x,\epsilon}$ when time large".  We certainly did not assume this, instead, this is exactly  the main goal of the  theorem. The statement   was  stated right after (3.29) in \cite{LR} and was  proved by the  contradiction argument  based on (3.22) in Proposition 3.4(\cite{LR}) ({ we copy this part below}).   In fact, this is  the main  technical issue  for Bott-Morse type gradient flows, though  the argument was  already standard for not-Morse type functions in gauge theory (cf. \cite{D},\cite{MMR}).
\v
We recite  our proof  in \cite{LR} from line 17  Page 175 to line 8  Page 176 below:
\v
\begin{mdframed}[skipabove=9pt,leftmargin=1cm,rightmargin=1cm]
{\em
We show that for any $C^{\infty}$-neighbourhood $U$ of $\{x_k(t+d),\;0\leq d\leq 1\}$ there is a $N>0$ such that if $s>N$ then $\tilde{u}(s,.)\in U.$ If not, we could find a neighbourhood $U\subset O$ and a subsequence of $s_i$ ( still denoted by $s_i$ ) and a sequence $s_i'$ such that
$$
\tilde{u}(s,.)\in O\;\; for\; s_i\leq s \leq s_i',\leqno(3.30)
$$
$$\tilde{u}(s_i',.)\notin U.\leqno(3.31)$$
 By Lemma (3.6) and by choosing a subsequence we may assume that
$$\tilde{u}(s_i',t)\rightarrow x'(k't)\;\;\;in \;C^{\infty}(S^1,\tilde{M})$$
for some $k'$-periodic solution $x'(t)\in O$. We may assume that $O$ is so small that there is no $k'$-periodic solution
in $O$ with $k'\ne k.$ Then, we have $k'=k$, $x'\in S_k$. We assume that $\tilde{E}(s)\ne 0$. From (3.29) we have
$$\int_{S^1}\tilde{d}(\tilde{u}(s,t), \tilde{u}(s_i,t))dt \leq \frac{1}{C}\left(\tilde{E}(s_i)\right)^{\tfrac{1}{2}},$$
where $\tilde{d}$ denotes the distance function defined by the metric $g_{\tilde{j}}$ on $Z$. Taking the limit $i\rightarrow \infty$, we get
$$\int_{S^1}\tilde{d}(x'(kt),x(kt))dt=0.$$
It follows that
$$x'(kt)=x(kt+\theta_0)$$
for some constant $\theta_0$. This contradicts (3.31). If there is some $s_0$ such that $\tilde{E}(s_0)=0$, then $|\Pi \tilde{u}_t|^2=|\Pi \tilde{u}_s|^2=0\;\;\forall s\geq s_0$. We still have a contradiction.
 }
\end{mdframed}

\v\v
 (3) Then  the argument for Theorem 3.7 can be applied to the case that $\mc A$  is Bott-Morse type, which  certainly {\em generalized } the result in \cite{HWZ}. If the case that the contact manifold
$\tilde M$ is a circle bundle of a line bundle, the similar theorem can be proved
in a rather easy way, for example, this point of view was employed  again  in \cite{CLSZ} and was  also used by  T\& Z. We remark  again, the way we adapt is suitable for more general cases for contact geometry.
\v
(4) It seems that T\& Z didn't know what we needed  later and didn't understand very well the standard elliptic estimates, even
didn't understand well  the results in  \cite{HWZ}.
In fact, let $C-D_1$ be a neighborhood of $\infty$ and $u:C-D_1\rightarrow \mathbb R\times \widetilde{M}$ be a J-holomorphic map, put $z=e^{s+2\pi i t},$  such that
\begin{enumerate}
\item[(i)] $E_{\phi}(u)\leq b;$
\item[(ii)] $|a(s,t)-ks-\ell_0|\rightarrow 0$,\;\;$|\theta(s,t)-kt-\theta_0|\rightarrow 0$,\;\;
$|y(s,t)|\rightarrow 0,$\;\;\;\;\; as \;\;$s\to +\infty;$
\item[(iii)] $u(s_0,t)$ lies in a compact set $K\subset R\times \tilde{M}$ where $s_0$ is a fixed point in $[1,+\infty)$.
\end{enumerate}
Then the constants $C_r$ in our paper depend only on $b$ and $K$.  This follows from the standard elliptic estimates and the estimates in [HWZ]. {  This uniform bound is certainly necessary for the gluing theory in \cite{LR}.}

\subsubsection{ }\label{sect_3.1.3}
On (LR6) and (LR8): T\& Z commented that  (1) the compactness argument in \cite{LR} is vague on the targets;
(2) a special type of maps, "contracted rubber map", is ignored.
\v
{\bf Response:}
\v
(1) From the context we actually assume that the target
of the sequence is $M^+$; the general cases can be dealt with similarly (for example, similar situations in Floer theory are well known);
\v
(2) "the contracted rubber map"  is {\em never} appeared for the cases studied  in \cite{LR}. As  we always considered  the map to $\mathbb R\times \widetilde M$, which can never be a contracted rubber map with only one puncture/node at one of the divisors, even in the compactification.

\subsubsection{ }\label{sect_3.1.4}
On (LR7): T\& Z questioned  (3) of Lemma 3.11 in \cite{LR}. They commented that "...in contrast to the setting in [H, HWZ1], the horizontal and vertical directions in the setting of [LR] are not tied together"(cf LR7)
"... the claim of [LR, Lemma 3.11(3)] in fact cannot be possibly true"(cf. Remark 4.5(\cite{TZ})).
\v

{\bf Response:}
\v
Our proof of (3) in  Lemma 3.11(\cite{LR}) is correct.   This is a standard geometry consequence from "no energy lost" in bubble tree argument and is proved by studying the energy on the "connecting neck" between two bubbles (for example, Lemma 4.5.1 and Page 57 in \cite{MS-1},  \S 6 in \cite{MMR}). It seems that T\&Z didn't understand  the standard  bubbling construction  in the literature very well.
Rather strangely,  T\&Z commented  that this is impossibly true. We sketch a simple  argument below.

By our construction (see (3.53) in \cite{LR})
$$v_i(r,t)=(b_i(r,t),\tilde{v}_i(r,t))= \left(a_i(log \delta_i+r,t)-a_i(log\delta_i,t_0),\tilde{u}_i(log \delta_i + r,t)\right).$$
This means that $\{u_i(-m,t),\;\;t\in S^1\}$ and $\{v_i(m,t),\;\;t\in S^1\}$ is connected by a tube. As $i\to \infty$, the connecting tube gives a gradient flow from one critical point
$\lim_{s\to-\infty}\tilde u(s,t)$ to another critical point $\lim_{r\to \infty}\tilde v(r,t)$. We have proved that the energy of the flow is 0, hence these two critical points are identical (in $W^2_r(S^1,\widetilde M)$).

\subsection{Gluing theory}\label{sect_3.2}
We answer (LR10-12) in this subsection.  {  Before we start with individual questions, we would like to comment on two general concerns raised by T\&Z.
The first concern  is the gluing theory for stable maps in rubber components. We understand, as we introduce an $\mathbb C^\ast$ action on moduli space on rubber components, it is fair to ask for a treatment accordingly. In fact, this is done by taking slices with respect to the action (cf. {\bf 2.} in Page 188 (\cite{LR})). Probably, T\&Z did not spot this?!  The second concern by T\&Z is how to generalize the gluing construction from one nodal case to general cases, which they think  might be nontrivial. We will  explain this point  in \S\ref{sect_3.2.3}.}

\subsubsection{ }\label{sect_3.2.1}
(LR10) comments on the gluing when the rubber components involved: (1) the gluing is not up to $\mathbb C^\ast$ action; (2) Remark 5.11(\cite{TZ}) concerns the target spaces with respect to the gluing parameters, in particular, when more than two gluing parameters are involved.
\v
{\bf Response:}
\v
(1) For this case, the slice is constructed  to deal with the  $\mathbb C^\ast$ quotient.
\v
In Page188 of \cite{LR}, from Line 26 to Line 31 we wrote:
\v
\begin{mdframed}[skipabove=9pt,leftmargin=1cm,rightmargin=1cm]
{\em {\bf 2:} $N=R\times \tilde{M}$. We must modulo the group $C^*$-action  generated by the $S^1$-action and the translation along $R$. We fix a point $y_0\in \sum$ different from the marked points and puncture points. Fix a local coordinate system $a,\theta, w$ on $R\times \tilde{M}$ such that $u(y_0)= (0,0,0)$. We use $\mc C'(\sum; u^*TN)=\{h\in \mc C(\sum;u^*TN)| h(y_0)=(0,0,*)\}$ instead of $\mc C(\sum;u^*TN)$, then the construction of the neighborhood $\mathfrak{u}_b$ of $b$ is the same as for $M^+$.}
\end{mdframed}
\v
This means that we constructed a slice, this equivalent to work on the quotient of the $C^*$-action.

\v
(2) The construction of gluing map itself  manifests  the parameter (the length of the connecting cylinder in the middle) for the target space explicitly. Though we did not spell out this parameter, it is fairly easy to get it. For example,  if only glue two
components $M^\pm$ with gluing parameter $(r,\theta)$, the length of connecting cylinder is $r$; if we glue three components $M^+,N$ and $M^-$(where $N$ is a rubber component between $M^\pm$) with two gluing parameters $(r^\pm,\theta^\pm)$, the length of connecting cylinder is
$r^++r^-$, and etc.

\subsubsection{(LR11) and (LR12)}\label{sect_3.2.2}
In (LR11,12), T\& Z complained  some routine issues in gluing theory such as (1) the injectivity and surjectivity  of the gluing maps were not explained, and (2) we applied  the implicit functional theorem without necessary bounds on the Taylor expansion of  $\bar\partial$.
\v
{\bf Response:}
\v
(1) { All arguments on injectivity and surjectivity are intensively developed by Taubes, Mrowka, Donaldson and etc in gauge theory and  already became  standard in  the gluing theory. In  the case of stable maps, this, for example, was  discussed by Fukaya-Ono ( see \cite{FO}, Chapter 3). Essentially, this does not require extra hard estimates.}   In this paper, we decided to mainly focus on the new issue that we concerned, for example, the estimates of $\frac{\p}{\p r}$ for gluing maps, {  which is certainly a {\em harder} issue. }
This refined estimate  was  certainly new at that time.
\begin{remark}
We would like to remark that in \cite{LR} the estimates for $\frac{\p}{\p r}$  is of order $r^{-2}$when $r\to \infty$, in fact, this can be achieved to be of exponential decay order $\exp(-Cr)$.
\end{remark}

\v
(2) This is a standard issue. Nevertheless, the necessary bounds for applying the implicit functional theorem were  given in Section 4.1. The $0$-th and $1$-th order were given by (4.16) and (4.26) in (\cite{LR}, P193, P195). The quadratic estimate was given by (4.3) in (\cite{LR}, P187).

\subsubsection{From one nodal case to general cases}\label{sect_3.2.3}
{  We first consider  the case with one gluing parameter, namely, we glue two components
of targets, say $M^\pm$, with gluing parameter $\rho=(r,\theta)$. Suppose we have a nodal surface consisting of $\Sigma^\pm$ with relative nodal points
$(p_1,\ldots, p_v)$ and a pair of relative maps $u^\pm:\Sigma^\pm\to M^\pm$. We want to glue a map $u_\rho: \Sigma_\rho\to M_\rho$. One of the main issue is to construct $\Sigma_\rho$.
The gluing parameter $\rho$ determines the gluing parameters $(\lambda_1(\rho),\ldots,\lambda_v(\rho))$ at nodal points. Here
$\lambda_i(\rho)$ is given by the local behaviour of $u$ at  $p_i$. We already gave this for $i=1$ (cf. (4.14),(4.15) in \cite{LR}), and of course, it can be generalized to multiple nodal points  in a straightforward way.

If the gluing parameters are more than one, this is already explained in \S\ref{sect_3.2.1}.
}

\subsection{Invariants}\label{sect_3.3}
We answer the comments (LR13,14,16) in this subsection.
\subsubsection{(LR13)}\label{sect_3.3.1}
In (LR13), T\& Z commented: VFC approach is based on a global regularization of the moduli space. The 3-4 pages dedicated to this could be avoided by using the local VFC approach of Fukaya-Ono or Li-Tian.
\v
{\bf Response:}
 \v
 Obviously,  the global regularization  if exists, would be much better than the  local regularization. The pay-off is 3-4 pages, while without the  local regularization, a further argument would be    of  more than 30-40 pages. So the logical comment is "we use global regularization to avoid the local VFC approach of  FO and LT".
\subsubsection{(LR14)}
In (LR14) T\& Z commented  that the construction for the rubber components needs to respect with the $\mathbb C^\ast$-action.
\v
{\bf Response:}
\v
 We always take slices for this group action.
\subsubsection{(LR16)}\label{sect_3.3.2}
In (LR16), T\& Z complained  that  it was not necessary to define the invariants using integrations
and, therefore, the estimates of $\frac{\p}{\p r}$ of gluing maps are not necessary .
\v
{\bf Response:} This is explained in \S\ref{sect_2.2}. We believe that T\&Z did not understand the main  idea and the contribution  in \cite{LR}.

\subsection{Others}\label{sect_3.4}
 Comments
(LR1), (LR2), (LR3), (LR9)  and (LR15)
 are on some miscellaneous points. We collect them at here.
\subsubsection{(LR1)}\label{sect_3.4.1}
{\bf Response to (LR1):}  in the formula of Theorem 5.8 (\cite{LR}), we agree with T\&Z
that we missed  an obvious term $\mathrm{Aut}(\mathbf k)$. This can be easily fixed by readers when they apply our formula.
\subsubsection{(LR2)}\label{sect_3.4.2}

{\bf Response to (LR2):} the way we formulate the relative stable maps followed \cite{R} for example. It might not be treated as  a standard way today though. The importance of this concept is that we first introduce the idea of "relative" and the  $\mathbb C^\ast$ action. T\&Z seem always try to ignore any significant point in \cite{LR} but focus on some non-essential points instead. Their attitude to the research paper under discussion  is certainly unprofessional, considering that they got some key ideas from \cite{LR}.

\subsubsection{(LR3)} \label{sect_3.4.3}
T\& Z wrote: "In addition to being imprecise, Definition 3.18 in \cite{LR} of the key notion of stable map of $X\cup_V Y$ ($\bar M^+\cup_D\bar M^-$ in the notation of \cite{LR}) is incorrect, as it seperates the rubber components into $X$ and $Y$-parts".
\v
{\bf Response:}
\v
By definition $M_{\infty}=M^+\bigcup M^-$ (corresponding to $\bar M^+\cup_D\bar M^-$ mentioned above).  There are two points of view of moduli spaces of maps in $M_\infty$.
\begin{enumerate}
\item[(\bf I)] Considering  the relative stable maps into $M_\infty$ with
matching condition at infinity, and then take the compactification. To be precise,
we have the moduli spaces $\mc M_{A^+}(M^+,g^+,T_{m^+})$, $\mc M_{A^-}(M^-,g^-,T_{m^-})$ and their compactification $\overline
{\mc M}_{A^+}(M^+,g^+,T_{m^+})$, $\overline{
\mc M}_{A^-}(M^-,g^-,T_{m^-})$. We can define the moduli space $\overline{\mc M}_A(M_{\infty},g,T_m)$ as a triple $(\Gamma^-,
\Gamma^+,\rho)$ where $\Gamma^{\pm}\in \overline{\mc M}_{A^\pm}(M^\pm,g^\pm,T_{m^\pm})$ and $\rho:\{p_1^+,...,p^+_v\}\rightarrow \{p_1^-,...,p^-_v\}$ is a one-to-one map satisfying the conditions described in Definition 3.18 in \cite{LR};
\item[(\bf J)] consider the relative stable maps into $M_r$ and take the limit
as $r\to \infty$.
\end{enumerate}

We know that T\&Z actually took  the approach ({\bf J}). The disadvantage of this is that this space is not virtually smooth. On the other hand, we know that ({\bf J})
fits into the degeneration framework well, for example, this has been  already used in
\cite{CLSZ} for orbifold cases.
The very important observation is, these two different choices have same top stratum, and thanks to our definition of invariants, they actually provide same invariants.
\begin{remark}\label{remark_3.4.1}
Moreover, the second condition (2) in Definition 3.18 makes it slightly different from the moduli spaces in the usual sense. Let us suppose that $v=1$. We have
 evaluation maps at $p_1^\pm$ to be $ev^\pm: \om_{A^\pm}\to D$. Define
$$
\mc M'=\overline{\mc M}_{A^+}(M^+,g^+,T_{m^+})\times_D\overline{\mc M}_{A^-}(M^-,g^-,T_{m^-})
$$
be the fiber product with respect to $ev^\pm$; however, Condition (2)
uses the "evaluation map" to space of closed orbits, it turns out that the
space $\overline{\mc M}_A(M_\infty,g,T_m)$ becomes the $k$ copies of
$\mc M'$, where $k$ is the multiplicity of the closed orbits at infinity.
\end{remark}

\subsubsection{(LR9)}\label{sect_3.4.4}
{\bf Response to (LR9):}
\v
 The Hausdorffness issue was  already standard in the subject even it may be nontrivial (see \cite{R}, Theorem 3.13 and \cite{FO}, Chapter 2). We do not think it is necessary to repeat standard materials, in particular for   papers  published in Invent. and Annuals.
\subsubsection{(LR15)}\label{sect_3.4.5}
{\bf Response to (LR15):}
\v
T\&Z asked the explanation of multiplicity $k$.  We agree that this needs more explanation. This $k$ factor is due to Remark \ref{remark_3.4.1} and it leads to a natural map $Q$ of degree $k$ used in P. 210(\cite{LR}).

\end {document}